\documentclass{amsart}
\usepackage{amsmath,amssymb}
\newtheorem{theorem}{Theorem}[section]
\newtheorem{lemma}[theorem]{Lemma}
\newtheorem{proposition}[theorem]{Proposition}
\newtheorem{corollary}[theorem]{Corollary}

\newtheorem{remark}[theorem]{Remark}

\begin{document}

\title{Injective modules and fp-injective modules over valuation rings}

\author{F. Couchot}

\begin{abstract} It is shown that each almost maximal valuation ring $R$,
such that every indecomposable injective $R$-module is countably
generated, satisfies the following condition (C): each fp-injective
$R$-module is locally injective. The converse holds if $R$ is a
domain. Moreover, it is proved that a valuation ring $R$ that
satisfies this condition (C) is almost maximal. The converse holds if
$\hbox{Spec}(R)$ is countable. When this last condition is satisfied it
is also proved that every ideal of $R$ is countably generated. New
criteria for a valuation ring to be almost maximal are given. They
generalize the criterion given by E. Matlis in the domain
case. Necessary and sufficient conditions for a valuation ring to be
an IF-ring are also given.
\end{abstract}
\maketitle

In the first part of this paper we study the valuation rings that
satisfy the following condition (C): every fp-injective module
is locally injective. In his paper \cite{fac}, Alberto Facchini constructs
an example of an almost maximal valuation domain satisfying (C) which
is not noetherian and gives a negative answer to the following
question asked in \cite{azu} by Goro Azumaya: if $R$ is a ring that
satisfies (C), is $R$ a left noetherian ring? From
\cite[Theorem 5]{fac} we easily deduce that a valuation domain $R$ satisfies
(C) if and only if $R$ is almost maximal and its classical field of
fractions is countably generated. In this case every indecomposable
injective $R$-module is countably generated. So, when an almost
maximal valuation ring $R,$ with eventually non-zero zerodivisors,
verifies this last condition, we prove that $R$ satisfies
(C). Conversely, every valuation ring that satisfies (C) is almost
maximal.

In the second part of this paper, we prove that every locally injective
module is a factor module of a direct sum of indecomposable injective
modules modulo a pure submodule. This result allows us to give equivalent
conditions for a valuation ring $R$ to be an IF-ring, i.e. a ring for
which every injective $R$-module is flat. It is proved that each
proper localization of $Q$, the classical ring of fractions of $R$, is
an IF-ring.

It is well known that a valuation domain $R$ is almost maximal if and
only if the injective dimension of the $R$-module $R$ is less or
equal to one. This result is due to E. Matlis. 
See \cite[Theorem 4]{mat}. In the third part, some generalizations of
this result are given. Moreover, when the subset $Z$ of zerodivisors
of an almost maximal valuation ring $R$ is nilpotent, we show that
every uniserial $R$-module is ``standard''(see \cite[p.141]{fuc}).

In the last part of this paper we determine some sufficient and necessary
conditions for every indecomposable injective module over a valuation
ring $R$ to be countably generated. In particular the following
condition is sufficient: $\mathrm{Spec}(R)$ is a countable set. Moreover, when this condition is satisfied, we
prove that every ideal of $R$ is countably generated and that every
finitely generated $R$-module is countably cogenerated.

In this paper all rings are associative and commutative with unity and
all mo\-dules are unital. An $R$-module $E$ is said to be \textit{locally
injective} (or finitely injective, or strongly absolutely pure, \cite{ram}) if every
homomorphism $A\rightarrow E$ extends to a homomorphism $B\rightarrow
E$ whenever $A$ is a finitely generated submodule of an arbitrary
$R$-module $B.$ As in \cite{facc} we say that $E$ is \textit{divisible} if, for every $r\in
R$ and $x\in E,$ $(0:r)\subseteq (0:x)$ implies that $x\in
rE$, and that $E$ is \textit{fp-injective}(or
absolutely pure) if $\mathrm{Ext}_R^1(F,E)=0,$ for every finitely
  presented $R$-module $F.$ A ring $R$ is called \textit{self
    fp-injective} if it is fp-injective as $R$-module. An exact sequence \ $0 \rightarrow F \rightarrow E \rightarrow G \rightarrow 0$ \ is \textit{pure}
if it remains exact when tensoring it with any $R$-module. In this case
we say that \ $F$ \ is a \textit{pure} submodule of $E$. Recall that a
module $E$ is fp-injective if and only if it is a pure submodule
of every overmodule (\cite[Proposition 2.6]{ste}). A module is said to be
\textit{uniserial} if its  submodules are linearly ordered by inclusion
and a ring $R$ is a \textit{valuation ring} if it is uniserial as
$R$-module. Recall that every finitely presented module over a
valuation ring is a finite direct sum of cyclic modules
\cite[Theorem 1]{war}. Consequently a module $E$ over a valuation ring
$R$ is fp-injective if and only if it is divisible. A valuation ring
$R$ is \textit{maximal} if every totally ordered family of cosets
$(a_i+L_i)_{i\in I}$ has a nonempty intersection and $R$ is
\textit{almost maximal} if the above condition holds whenever $\cap_{i\in
  I}L_i\ne 0.$ 

We denote $\hbox{p.d.}_R\,M$ (resp. $\hbox{i.d.}_R\,M$) the projective
(resp. injective) dimension of $M,$ where $M$ is a
module over a ring $R,$ $E_R(M)$ the injective hull of $M,$
$\hbox{Spec}(R)$ the space of prime ideals of $R,$ and for every ideal
$A$ of $R,$ $V(A)=\{I\in\hbox{Spec}(R)\mid A\subseteq I\}$
  and $D(A)=\hbox{Spec}(R)\setminus V(A).$

When $R$ is a valuation ring, we denote by $P$ its maximal ideal, $Z$ its
subset of zerodivisors and $Q$ its classical ring of fractions. Then $Z$
is a prime ideal and $Q=R_Z.$

\section{Valuation rings whose fp-injective modules are locally injective}
\label{S:loin}

From \cite[Theorem 5]{fac} we easily deduce the following theorem:
\begin{theorem} Let $R$ be a valuation domain. The following assertions
  are equivalent:
\begin{enumerate}
\item Every fp-injective module is locally injective.
\item $R$ is almost maximal and $\hbox{p.d.}_R\,Q\leq 1.$
\item $R$ is almost maximal and $Q$ is countably generated over $R.$
\item $R$ is almost maximal and every indecomposable injective module
  is countably generated.
\end{enumerate}
\end{theorem}
\textbf{Proof.} $(1)\Rightarrow (2).$ Since $Q/R$ is fp-injective, it is a
locally injective module. By \cite[Corollary 3.4]{ram}, $Q/R$ is
injective and consequently $R$ is almost maximal by 
\cite[Theorem 4]{mat}. From \cite[Theorem 5]{fac}, it follows that
$\hbox{p.d.}_R\,Q\leq 1.$ 

$(2)\Rightarrow (1)$ is proved in \cite[Theorem 5]{fac}.

$(2)\Leftrightarrow (3)$. See \cite[Theorem 2.4, p.76]{fuc}.

$(3)\Leftrightarrow (4).$ By \cite[Theorem 4]{mat} 
$E_R(R/A)\simeq Q/A$ for every proper ideal $A.$ Consequently,
if $Q$ is countably generated, every indecomposable injective module
is also countably generated.\qed
 
\bigskip
If $R$ is not a domain then the implication $(4)\Rightarrow (1)$
holds. The following lemma is needed to prove this implication and
will be useful in the sequel too.
\begin{lemma}
\label{L:A}
Let $R$ be a valuation ring, $M$ an $R$-module,
$r\in R$ and $y\in M$ such that $ry\not= 0.$ Then:
\begin{enumerate}
\item $(0:y)=r(0:ry).$
\item If $(0:y)\not= 0$ then $(0:y)$ is finitely
  generated if and only if $(0:ry)$ is finitely generated.
\end{enumerate}
\end{lemma}
\textbf{Proof.} Clearly $r(0:ry)\subseteq (0:y).$ Let
$a\in (0:y).$ Since $ry\not= 0$, $(0:y)\subset rR.$ There
exists $t\in R$ such that $a=rt$ and we easily check that $t\in
(0:ry).$ The second assertion is an immediate consequence of
the first. \qed
 
\begin{theorem} Let $R$ be an almost maximal valuation
ring. Assume that eve\-ry indecomposable injective $R$-module is
countably generated.
Then every fp-injective $R$-module is locally injective.
\end{theorem}
\textbf{Proof.} Let $F$ be a non-zero fp-injective module. We must prove
that $F$ contains an injective hull of each of its finitely generated
submodules by \cite[Proposition 3.3]{ram}. Let $M$ be a finitely
generated submodule of $F.$ By \cite[Theorem]{gil} $M$ is a
finite direct sum of cyclic submodules. Consequently, we may assume
that $M$ is cyclic, generated by $x$. Let $E$ be an injective hull of
$M$ and $\{x_n\mid n\in \mathbb{N}\}$ a spanning set of $E.$ By \cite[Theorem]{gil} $E$ is a uniserial module. Hence, for every
integer $n,$ there exists $c_n\in R$ such that $x_n=c_nx_{n+1}.$ We
may suppose that $x=x_0.$ By induction on $n$ we prove that there
exists a sequence $(y_n)_{n\in\mathbb{N}}$ of elements of $F$ such that
$y_0=x,$ $(0:x_n)=(0:y_n)$ and $y_n=c_ny_{n+1}.$
Since $x_n=c_nx_{n+1},$ $(0:c_n)\subseteq
(0:x_n)=(0:y_n).$ Since $F$ is fp-injective, there
exists $y_{n+1}\in F,$ such that $y_n=c_ny_{n+1}.$ We easily deduce 
  from Lemma \ref{L:A} that $(0:x_{n+1})=(0:y_{n+1}).$ Now,
  the submodule of $F$ generated by $\{y_n\mid n\in \mathbb{N}\}$ is
  isomorphic to $E.$ \qed
 
\bigskip
We don't know if the converse of this theorem holds when $R$ is not a
domain. However, for every valuation ring $R$, condition (C)
implies that $R$ is almost maximal. 
 
\begin{theorem}
\label{T:A}
Let $R$ be a valuation ring. 
If every fp-injective $R$-module is locally injective then $R$ is
almost maximal.
\end{theorem}

Some preliminary results are needed to prove this theorem. The
following Lemma will often be used in the sequel. This lemma is
similar to \cite[Lemma~II.2.1]{fuc}.
\begin{lemma}
\label{L:B} 
Let $R$ be a local commutative ring, $P$ its
maximal ideal, $U$ a uniserial $R$-module, $r\in R,$ and $x,y\in U$
such that $rx=ry\ne 0.$ Then $Rx=Ry.$
\end{lemma}
\textbf{Proof.} We may assume that $x=ty$ for some $t\in R.$ It follows
that $(1-t)ry=0.$ Since $ry\ne 0$ we deduce that $t$ is a unit. \qed
 
\begin{proposition} 
\label{P:A}
Let $U$ be a uniform fp-injective module over a valuation ring $R$. Suppose
   there exists a nonzero element $x$ of $U$ such that $Z=(0:x)$. Then:
\begin{enumerate}
\item $U$ is a $Q$-module.
\item For every proper $R$-submodule $A$ of $Q$, $U/Ax$ is faithful and
  fp-injective.
\end{enumerate}
\end{proposition}
\textbf{Proof.} (1) For every $0\ne y\in U$, $(0:y)=sZ$ or
$(0:y)=(Z:s)=Z$ (see \cite{nis}). Hence
$(0:y)\subseteq Z$. If $s\in R\setminus Z$ then
the multiplication by $s$ in $U$ is injective. Since
$U$ is fp-injective this multiplication is bijective.

(2) If $R\subseteq A$ there exists $s\in R\setminus Z$ such that
$sA\subset R$ and there exists $y\in U$ such that $x=sy$. Then
$Ax=Asy$ and $(0:y)=Z$. Consequently we may assume that
$A\subset R$, after eventually replacing $A$ with $As$ and $x$ with $y$. Let $t\in R$. Since $(0:t)\subseteq Z$ there exists
$z\in U$ such that $x=tz$. Therefore $0\ne x+Ax=t(z+Ax)$ whence
$U/Ax$ is faithful. Let $t\in R$ and $y\in U$ such that
$(0:t)\subseteq (0:y+Ax)$. Therefore $(0:t)y\subseteq Ax\subset
Qx$. It is easy to check that $(0:t)$ is an ideal of $Q$. Since $Qx$
is the nonzero minimal $Q$-submodule of $U$ we get that $(0:t)\subseteq
(0:y)$. Since $U$ is fp-injective we conclude that $U/Ax$ is
fp-injective too. \qed
 
\bigskip
Now, we prove Theorem~\ref{T:A}.

\textbf{Proof of Theorem~\ref{T:A}.} If $Z=P$ then $R$ is self
  fp-injective by \cite[Lemma 3]{gil}. It follows that $R$ is self
  injective by \cite[Corollary 3.4]{ram} and that $R$ is maximal by
  \cite[Theorem 2.3]{kla}.

Now we assume that $Z\ne P.$ In the
same way we prove that $Q$ is maximal. From \cite[Theorem]{gil} it follows that $E_Q(Q/Z)\simeq E_R(R/Z)$ is uniserial over
$Q$ and $R.$ Let $H=E_R(R/Z)$ and $x\in H$ such that $Z=(0:x)$.
By Proposition~\ref{P:A} $H/Px$ is fp-injective. This module is injective by
\cite[Corollary 3.4]{ram}. Hence $E(R/P)\simeq H/Px$ is uniserial. By
\cite[Theorem]{gil} $R$ is almost maximal.  \qed
  
\bigskip
From Proposition~\ref{P:A} we  easily deduce the following corollary which
generalizes the second part of \cite[Theorem 4]{mat}.
\begin{corollary}
\label{C:A} 
Let $R$ be an almost maximal valuation ring,
$H=E(R/Z)$ and $x\in H$ such that $Z=(0:x)$.
For every proper and faithful ideal $A$ of $R$, $H/Ax\simeq E(R/A)$.
\end{corollary}
\textbf{Proof.} By \cite[Theorem]{gil} $E(R/A)$ is uniserial. It
  follows that its proper submodules are not faithful. We conclude by
  Proposition~\ref{P:A}. \qed
 
\section{Valuation rings that are IF-rings}
\label{S:ifr}

We begin this section with some results on indecomposable injective
modules over a valuation ring. In the sequel, if $R$ is a valuation
ring, let $E=E(R)$, $H=E(R/Z)$ and $F=E(R/Rr)$ for every $r\in
P,$ $r\not= 0.$ Recall that, if $r$ and $s$ are nonzero elements of
$P,$ then $E(R/Rr)\simeq E(R/Rs)$, (see \cite{nis}).
\begin{proposition} 
\label{P:fla}
The following statements hold for a valuation ring $R$.
\begin{enumerate}
\item The modules $E$ and $H$ are flat.
\item The modules $E$ and $H$ are isomorphic if and only if $Z$ is not
faithful.
\end{enumerate}
\end{proposition}
\textbf{Proof.} First we assume that $Z=P,$ whence
$R$ is fp-injective. Let $x\in E,$ $x\ne 0,$ and $r\in R$ such that
$rx=0.$ There exists $a\in R$ such that $ax\in R$ and $ax\not= 0.$
Then $(0:a)\subseteq (0:ax)$, so that there exists $d\in R$ such that
$ax=ad.$ By Lemma \ref{L:A} $(0:d)=(0:x),$ whence
there exists $y\in E$ such that $x=dy.$ We deduce that $r\otimes
x=rd\otimes y=0.$ Hence $E$ is flat. Now if $Z\ne P,$ then
$E\simeq E_Q(Q).$ Consequently $E$ is flat over $Q$ and $R.$

Since $Q$ is self-fp-injective, $E_Q(Q/Z)\simeq H$ is flat by \cite[Theorem
2.8]{cou}.

If $Z$ is not faithful there exists $a\in Z$ such that $Z=(0:a)$. It follows
that $H\simeq E(Ra)=E$. \qed 
 
\bigskip
We state that $E$ and $F$ are generators of the category of locally
injective $R$-modules. More precisely:
\begin{proposition} 
\label{P:gen}
Let $R$ be a valuation ring and $G$ a
locally injective module. Then there exists a pure exact sequence:
$0\rightarrow K\rightarrow I\rightarrow G\rightarrow 0$, such
that $I$ is a direct sum of submodules isomorphic to $E$ or $F.$ 
\end{proposition}
\textbf{Proof.} There exist a set $\Lambda$ and an epimorphism
$\varphi:L=\oplus_{\lambda\in\Lambda}R_{\lambda}\rightarrow G$, where
$R_{\lambda}=R,\forall \lambda\in\Lambda.$ Let
$u_{\mu}:R_{\mu}\rightarrow L$ the canonical monomorphism. For every
$\mu\in\Lambda,$ $\varphi\circ u_{\mu}$ can be extended to
$\psi_{\mu}:E_{\mu}\rightarrow G,$ where
$E_{\mu}=E,\forall\mu\in\Lambda.$ We denote
$\psi:\oplus_{\mu\in\Lambda}E_{\mu}\rightarrow G,$ the epimorphism defined by
the family $(\psi_{\mu})_{\mu\in\Lambda}.$ We put
$\Delta=\mathrm{hom}_R(F,G)$ and $\rho:F^{(\Delta)}\rightarrow G$ the morphism
defined by the elements of $\Delta.$ Thus $\psi$ and $\rho$ induce an
epimorphism $\phi:I=E^{(\Lambda)}\oplus F^{(\Delta)}\rightarrow G.$ Since, for
every $r\in P,r\ne 0,$ each morphism $g:R/Rr\rightarrow G$ can be extended to
$F\rightarrow G,$ we deduce that $K=\ker \phi$ is a pure
submodule of $I.$ \qed
 
\bigskip
Recall that a ring $R$ is \textit{coherent} if every finitely generated
ideal of $R$ is finitely presented. As in \cite{col} we say that $R$ is
an \textit{IF-ring} if every injective $R$-module is flat. From 
Propositions \ref{P:fla} and \ref{P:gen} we deduce necessary and
sufficient conditions for a valuation ring to be an IF-ring.
\begin{theorem} 
\label{T:ifr}
Let $R$ be a valuation ring which is not a
field. Then the following assertions are equivalent:
\begin{enumerate}
\item $R$ is coherent and self-fp-injective
\item $R$ is an IF-ring
\item $F$ is flat
\item $F\simeq E$
\item $P$ is not a flat $R$-module
\item There exists $r\in R,r\ne 0,$ such that $(0:r)$ is a nonzero
  principal ideal.
\end{enumerate}
\end{theorem}
\textbf{Proof.}(1)$\Rightarrow$(4). By \cite[Corollary 3]{col}, for every
$r\in P,r\ne 0,$ there exists $t\in P,t\ne 0,$ such that $(0:t)=Rr.$
Hence $R/Rr\simeq Rt\subseteq R\subseteq E.$ We deduce that $F\simeq
E.$

(4)$\Rightarrow$(3) follows from Proposition \ref{P:fla}.

(3)$\Rightarrow$(2) If $G$ is an injective module then by Proposition \ref{P:gen}
   there exists a pure exact sequence $0\rightarrow K\rightarrow
   I\rightarrow G\rightarrow 0$ where $I$ is a direct sum of
   submodules isomorphic to $E$ or $F$. By Propositions \ref{P:fla} $I$ is
   flat whence $G$ is flat too.

(2)$\Rightarrow$(1). See \cite[Theorem 2]{col}.

(1)$\Rightarrow$(6) is an immediate consequence of \cite[Corollary 3]{col}.

(6)$\Rightarrow$(5) We denote $(0:r)=Rt.$ If $r\otimes t=0$ in
$Rr\otimes P$ then, by \cite[Proposition~13, p. 42]{bou}, there exist $s$
and $d$ in $P$ such that $t=ds$ and $rd=0.$ Thus $d\in (0:r)$ and
$d\notin Rt.$ Whence a contradiction. Consequently $P$ is not flat.

(5)$\Rightarrow$(1). If $Z\ne P,$ then $P=\cup_{r\notin Z}Rr$,
whence $P$ is flat. Hence $Z=P.$ If $R$ is not coherent, there exists
$r\in P$ such that $(0:r)$ is not finitely generated. By Lemma \ref{L:A}
$(0:s)$ is not finitely generated for each $s\in P,s\ne
0.$ Consequently, if $st=0$, there exist $p\in P$ and $a\in (0:s)$
such that $t=ap.$ It follows that $s\otimes t=sa\otimes p=0$ in
$Rs\otimes P.$ Whence $P$ is a flat module. We get a contradiction.
\qed
 
\bigskip
The following theorem allows us to give examples of valuation rings that are
IF-rings.
\begin{theorem} 
\label{T:ifl}
The following statements hold for a valuation ring $R$:
\begin{enumerate}
\item For every $0\not= r\in P$, $R/Rr$ is an IF-ring.
\item For every prime ideal $J\subset Z$, $R_J$ is an IF-ring.
\end{enumerate}
\end{theorem}
\textbf{Proof.} (1) For every $a\in P\setminus Rr$ there exists $b\in
P\setminus Rr$ such that $r=ab$. We easily deduce that $(Rr:a)=Rb$
whence $R/Rr$ is an IF-ring by Theorem~\ref{T:ifr}.

(2) The inclusion $J\subset Z$ implies that there exist $s\in
    Z\setminus J$ and $0\ne r\in J$ such that $sr=0$. If we set $R'=R/Rr$
    then $R_J\simeq R'_J$. From the first part and 
    \cite[Proposition 1.2]{cou} it follows that $R_J$ is an IF-ring.
\qed
 
\bigskip
The two following lemmas are needed to prove the important
Proposition \ref{P:main}.
\begin{lemma} 
\label{L:E}
The following statements hold for a valuation ring $R$:
\begin{enumerate} 
\item If $Z\not=P$ then $E=PE$.
\item If $Z=P$ then $E=R+PE$ and $E/PE\simeq R/P$.
\end{enumerate}
\end{lemma}
\textbf{Proof.} (1) If $p\in P\setminus Z$ then $E=pE$.

(2) For every $x\in PE$, $(0:x)\ne 0$ whence $1\notin
PE$. Let $x\in E\setminus R$. There exists $r\in R$ such that
$0\not=rx\in R$. Since $R$ is self-fp-injective there exists $d\in R$
such that $rd=rx$. By Lemma~\ref{L:A} $(0:d)=(0:x)$. We deduce that
$x=dy$ for some $y\in E$. Then $x\in PE$ if $d\in P$. If $d$ is a
unit, in the same way we find $t,c\in R$ and $z\in E$ such that
$tc=t(x-d)\not=0$ and $x-d=cz$. Since $r\in (0:x-d)=(0:c)$ then
$c\in P$ and $x\in R+PE$. \qed
 
\begin{lemma} 
\label{L:unif}
Let $R$ be a valuation ring and $U$ a uniform
$R$-module. If $x,y\in U$, $x\notin Ry$ and $y\notin Rx$, then $Rx\cap
Ry$ is not finitely generated.
\end{lemma}
\textbf{Proof.} Suppose that $Rx\cap Ry=Rz$. We may assume that there exist
  $t\in P$ and $d\in R$ such that $z=ty=tdx$. It is easy to check that
$(Rx:y-dx)=(Rx:y)=(Rz:y)=Rt\subseteq (0:y-dx)$. It follows that
$Rx\cap R(y-dx)=0$. This contradicts that $U$ is uniform.
\qed
 
\begin{proposition} 
\label{P:main}
Let $R$ be a valuation ring which is not a
field. Apply the functor $\mathrm{Hom}_R(-,E(R/P))$ to the canonical
exact sequence 
\[(S): 0\rightarrow P\rightarrow R\rightarrow
R/P\rightarrow 0.\ Then:\]
\begin{enumerate}
\item If $R$ is not an IF-ring one gets an exact sequence
\[(S_1): 0\rightarrow R/P\rightarrow E(R/P)\rightarrow F\rightarrow
0,\ with\ F\simeq \mathrm{Hom}_R(P,E(R/P)),\]
\item If $R$ is an IF-ring one gets an exact sequence
\[(S_2): 0\rightarrow R/P\rightarrow E(R/P)\rightarrow F\rightarrow
R/P\rightarrow 0,\ with\ PF\simeq \mathrm{Hom}_R(P,E(R/P)).\]
\end{enumerate}
\end{proposition}
\textbf{Proof.} (1) (S) induces the following exact sequence:

$0\rightarrow R/P\rightarrow E(R/P)\rightarrow \mathrm{Hom}_R(P,E(R/P))\rightarrow
0$. By Theorem \ref{T:ifr} $P$ is flat whence $\mathrm{Hom}_R(P,E(R/P))$ is
injective. Let $f$ and $g$ be two nonzero elements of
$\mathrm{Hom}_R(P,E(R/P))$. There exist $x$ and $y$ in $E(R/P)$ such
that $f(p)=px$ and $g(p)=py$ for each $p\in P$. Let $Rv$ be the minimal
nonzero submodule of $E(R/P)$. By Lemma \ref{L:unif} there exists $z\in (Rx\cap
Ry)\setminus Rv$. Then the map $h$ defined by $h(p)=pz$ for each
$p\in P$ is nonzero and belongs
to $Rf\cap Rg$. Thus $\mathrm{Hom}_R(P,E(R/P))$ is uniform. Now let
$a\in R$ such that $af=0$. It follows that
$Pa\subseteq (0:x)=Pb$ for some $b\in R$. We deduce that
$(0:f)=Rb$. Hence $F\simeq \mathrm{Hom}_R(P,E(R/P))$.

(2) First we suppose that $P$ is not finitely generated. From the first part of the proof it follows that
    $\mathrm{Hom}_R(P,E(R/P))\subseteq F$. We use the same
    notations as in (1). We have $(0:f)=Rb$ and
    there exists $c\in P$ such that $(0:c)=Rb$.  Consequently $f\in
    cF\subseteq PF$. Conversely let $y\in PF$ and $b\in P$ such that
    $(0:y)=Rb$. Since $R'=R/Pb$ is not an IF-ring it follows from
    the first part that $\mathrm{Hom}_R(P/bP,E(R/P))\simeq
    \{x\in F\mid bP\subseteq (0:x)\}\simeq E_{R'}(R'/rR')$ where
    $0\ne r\in P/bP$. Hence $\mathrm{Hom}_R(P,E(R/P))=PF$. We
    deduce the result from Theorem \ref{T:ifr} and Lemma \ref{L:E}.

If $P=pR$ then $E(R/P)\simeq E\simeq F$. Then multiplication by $p$
induces the exact sequence $(S_2)$. \qed
 
\begin{remark} 
\label{R:fac}
\textnormal{In \cite[Theorem 5.7]{facc} A. Facchini considered
    indecomposable pure-injective modules over a valuation ring
    $R$. He proved that $\mathrm{Hom}_R(W,G)$ is indecomposable for
    every indecomposable injective $R$-module $G$, where $W$ is a faithful
    uniserial module such that $(0:x)$ is a nonzero principal ideal
    for each $x\in W$. This result implies that
    $\mathrm{Hom}_R(P,E(R/P))$ is indecomposable when $R$ is an
    IF-ring and $P$ is faithful.}
\end{remark}

From Proposition~\ref{P:main} we deduce a sufficient and necessary condition
for a valuation ring to be an IF-ring. As in \cite{ste}, the
\textit{fp-injective dimension} of an $R$-module $M$
($\mathrm{fp-i.d.}_R\,M$) is the smallest integer $n\geq 0$ such that
$\hbox{Ext}^{n+1}_R(N,M)=0$ for every finitely presented $R$-module
$N$.
\begin{corollary} 
\label{C:ifr}
Let $R$ be a valuation ring. Then the
following assertions are equivalent:
\begin{enumerate}
\item $R$ is not an IF-ring. 
\item $\mathrm{i.d.}_R\,R/P=1$.
\item $\mathrm{fp-i.d.}_R\,R/Z=1$.
\end{enumerate}
\end{corollary}
\textbf{Proof.} $(1)\Leftrightarrow (2)$. It is an immediate consequence
of Proposition \ref{P:main}.

$(1)\Leftrightarrow (3)$. If $R$ is an IF-ring then $Z=P$ and
$\mathrm{fp-i.d.}_R\,R/Z>1$ by Proposition~\ref{P:main}. Assume that $R$
is not an IF-ring and $Z\ne P$. Let $x\in H$ such
that $Z=(0:x)$. By Proposition \ref{P:A} $H/Rx$ is fp-injective. It
follows that $\mathrm{fp-i.d.}_R\,R/Z=1$.  \qed
 
\section{Injective modules and uniserial modules}
\label{S:uni}

Proposition \ref{P:main} allows us to give generalizations of
well known results in the domain case. This is a first generalization
of the first part of \cite[Theorem 4]{mat}. 
\begin{theorem} 
\label{T:mat}
Let $R$ be a valuation ring. Then $R$ is
almost maximal if and only if $F$ is uniserial.
\end{theorem}
\textbf{Proof.} By \cite[Theorem]{gil} $F$ is uniserial if $R$
is almost maximal.  Conversely if $F$ is uniserial, by using the exact
sequence $(S_1)$ or $(S_2)$ of Proposition \ref{P:main}, it is easy to prove that $E(R/P)$ is
uniserial. We conclude by using \cite[Theorem]{gil}. 
\qed
 
\bigskip
Now we shall prove the existence of uniserial fp-injective
modules. This is an immediate consequence of \cite[Theorem]{gil} when $R$ is an almost maximal valuation ring. The
following proposition will be useful for this. 
\begin{proposition} 
\label{P:cof}
Let $R$ be a commutative local ring, $P$
its maximal ideal, $U$ a uniserial
module with a nonzero minimal submodule $S$. 
Then $U$ is a divisible module if and only if it is faithful.
\end{proposition}
\textbf{Proof.} First we suppose that $U$ is faithful. Let $0\ne s\in P$
    and $0\ne y\in U$ such that $(0:s)\subseteq
    (0:y)$. There exists $t\in R$ such that $x=ty$ where $x$
    generates $S$. Thus $t\notin (0:s)$ and consequently $stU$ is a
    nonzero submodule of $U$. It follows that there exists $z\in U$
    such that $x=stz$. By Lemma \ref{L:B} $y\in sU$.
Conversely let $0\ne s\in P$. Then $(0:s)\subseteq P$ implies that
$x\in sU$. We conclude that $U$ is faithful. \qed
 
\begin{proposition} 
\label{P:U}
Let $R$ be a valuation ring such that
$Z=P$. Assume that $R$ is coherent, or $(0:P)\not=0$, or $0$
is a countable intersection of nonzero ideals.
Then there exist two uniserial fp-injective modules $U$ and $V$ such
that $E(U)\simeq E(R/P)$ and $E(V)\simeq F$. When $P$ is principal
then $U\simeq V\simeq R$.
\end{proposition}
\textbf{Proof.} If $P$ is principal then $R$ is an IF-ring and it is
obvious that $U\simeq V\simeq R$. 

Now suppose that $R$ is an IF-ring
and $P$ is not finitely generated. Consequently $P$ is faithful. Let
$\phi:E(R/P)\rightarrow F$ be the
homomorphism defined in Proposition \ref{P:main}. Thus $F\simeq E$ and
$\mathrm{Im }\phi=PF$. It is obvious that $PF\cap R=P$. We put $V=R$ and
$U=\phi^{-1}(P)$. Since $P$ is faithful, $U$ is faithful too. It is
easy to prove that $U$ is uniserial. Then $U$ is fp-injective by
Proposition \ref{P:cof}. 

Now we assume that $P$ is not faithful and not finitely
generated. Then $R$ is not an IF-ring. By Corollary \ref{C:ifr}
$\mathrm{i.d.}_R\,(R/P)=1$. It follows that $R/(0:P)$ is
fp-injective. In this case we put $U=R$ and $V=R/(0:P)$.

Now we suppose that $0$ is a countable intersection of nonzero
ideals. We may assume that $R$ is not coherent and
$P$ is faithful. By \cite[Theorem 5.5]{sho} there exists a faithful
uniserial $R$-module $U$ such that $E(U)\simeq E(R/P)$. By Proposition
\ref{P:cof} $U$ is fp-injective. Let $u\in U$ such that $(0:u)=P$. Since
$R$ is not an IF-ring, then by using Corollary \ref{C:ifr} it is easy to
prove that $U/Ru$ is fp-injective. We put $V=U/Ru$. \qed
 
\begin{remark} \textnormal{ By \cite[Theorem 5.4]{facc} it is obvious
    that every faithful indecomposable pure-injective $R$-module is
    injective if $R$ is a valuation ring such that $(0:P)\ne
    0$. In this case the module $W$ in remark~\ref{R:fac} doesn't
    exist. By \cite[Theorem 5.7]{facc} and Proposition~\ref{P:main},
    $PF$ is the only faithful indecomposable
    pure-injective $R$-module which is not injective, when $R$ is an
    IF-ring and $(0:P)=0$.}
\end{remark}

As in \cite[p.15]{fuc}, for every proper ideal $A$ of a valuation ring
$R$ we put\\ 
$A^{\#}=\{s\in R\mid (A:s)\not=A\}.$ Then $A^{\#}/A$ is the set of
zerodivisors of $R/A$ whence $A^{\#}$ is a prime ideal. In particular
$0^{\#}=Z$. 
\begin{lemma} 
\label{L:die}
Let $R$ be a valuation ring, $A$ a proper ideal
of $R$ and $t\in R\setminus A$. Then $A^{\#}=(A:t)^{\#}$.
\end{lemma}
\textbf{Proof.} Let $a\in (A:t)^{\#}$. If $a\in (A:t)$ then $a\in
  A^{\#}$. If $a\notin (A:t)$ there exists $c\notin (A:t)$ such that
  $ac\in (A:t)$. If follows that $act\in A$ and $ct\notin A$ whence $a\in
A^{\#}$. Conversely let $a\in A^{\#}$. There exists $c\notin A$ such
that $ac\in A$. If $a\in (A:t)$ then $a\in (A:t)^{\#}$. If $a\notin
(A:t)$ then $at\notin A$. Since $ac\in A$ it follows that $c=bt$ for
some $b\in P$. Since $c\notin A$ it follows that $b\notin (A:t)$. From
$abt\in A$ we successively deduce that $ab\in (A:t)$ and $a\in
(A:t)^{\#}$. \qed
 
\bigskip
If $J$ is a prime ideal contained in $Z$, we put $ke(J)$ the kernel of
the natural map: $R\rightarrow R_J$.
\begin{corollary} 
\label{C:U}
Let $R$ be a valuation ring. Then:
\begin{enumerate}
\item For every prime ideal $J\subset Z$ there exist two uniserial
 fp-injective modules $U_{(J)}$ and $V_{(J)}$ such that
   $E(U_{(J)})\simeq E(R/J))$ and $E(V_{(J)})\simeq E(R_J/rR_J)$, where
   $r\in J\setminus ke(J)$.
\item If $Q$ is coherent, or $Z$ is not faithful, or $0$ is a
  countable intersection of nonzero ideals, there exist two uniserial
  fp-injective mo\-dules $U_{(Z)}$ and $V_{(Z)}$ such that
  $E(U_{(Z)})\simeq H$ and $E(V_{(Z)})\simeq E(Q/rQ)$, where
  $0\not=r\in Z$.  
\item If $Q$ is coherent, or $Z$ is not faithful, or $0$
is a countable intersection of nonzero ideals, then for every proper ideal $A$
  such that $Z\subset A^{\#}$ there exists a faithful uniserial
  fp-injective module $U_{(A)}$ such that $E(U_{(A)})\simeq E(R/A)$.
\end{enumerate}
\end{corollary}
\textbf{Proof.} (1) is a consequence of Theorem~\ref{T:ifl} and
  Proposition~\ref{P:U}.

(2) follows from the above proposition.

(3) is a consequence of (2) and Proposition~\ref{P:A}. More precisely, since
    $Z\subset A^{\#}$ we may assume that $A$ is faithful,
    eventually after replacing $A$ with $(A:a)$ for some $a\in
    A^{\#}$. Then we put $U_{(A)}=U_{(Z)}/Au$ where $u\in
    U_{(Z)}$ and $(0:u)=Z$.
\qed
 
\bigskip
From Theorem~\ref{T:mat}, Corollary~\ref{C:A}, Corollary~\ref{C:U} and
Corollary~\ref{C:ifr} we deduce another genera\-lization of
\cite[Theorem 4]{mat}.
\begin{theorem} 
\label{T:matl}
Let $R$ be a valuation ring. Suppose that $Q$
    is coherent or ma\-ximal, or $Z$ is not faithful, or $0$ is a countable
    intersection of nonzero ideals. Let $U_{(Z)}$ be the fp-injective
    uniserial module defined in Corollary~\ref{C:U} and $u\in U_{(Z)}$ such
    that $Z=(0:u)$. Then: 
\begin{enumerate}
\item If $R$ is not an IF-ring then $R$ is almost maximal if and only
  if $U_{(Z)}/Ru$ is injective. 
\item If $R$ is almost maximal then for every proper and faithful
  ideal $A$ of $R$, $U_{(Z)}/Au\simeq E(R/A)$.
\end{enumerate}
\end{theorem}
\textbf{Proof.} (1) If $U_{(Z)}/Ru$ is injective then $F\simeq U_{(Z)}/Ru$ is
  uniserial. By Theorem~\ref{T:mat} $R$ is almost maximal. Conversely,
  $U_{(Z)}/Ru$ is a fp-injective submodule of $F$ by
  Corollary~\ref{C:ifr} and $F$ is uniserial by
  \cite[Theorem]{gil}. Let $0\ne x\in F$. There exists $a\in R$ such
  that $0\ne ax\in U_{(Z)}/Ru$. It follows that $\exists y\in
  U_{(Z)}/Ru$ such that $ax=ay$. By using Lemma~\ref{L:B} we deduce
  that $x\in U_{(Z)}/Ru$. Hence $U_{(Z)}/Ru$ is injective.
  
(2) is an immediate consequence of Corollary~\ref{C:A}  \qed
 
\bigskip
Let us observe that $U_{(Z)}=Q$ when $Z$ is not faithful and
consequently we have a generalization of \cite[Theorem 4]{mat}.

Now this is a generalization of \cite[Theorem VII.1.4]{fuc}.
\begin{theorem} Let $R$ be an almost maximal valuation ring and
   suppose that $Z$ is nilpotent. Then:
\begin{enumerate} 
\item Every indecomposable injective $R$-module is a faithful factor
  of $Q$.
\item An $R$-module $U$ is uniserial if and only if it is of the form
  $U\simeq J/I$ where $I\subset J$ are $R$-submodules of $Q$. 
\item If $I\subset J$ and $I'\subset J'$ are $R$-submodules of $Q$
  then $J/I\simeq J' /I'$ if and only if $I=(I':q)$ and $J=(J':q)$ for some
  $0\ne q\in Q$.
\end{enumerate}
\end{theorem}
\textbf{Proof.} (1) $Q$ is an artinian ring and $Z$ is its unique prime
ideal. If $A$ is an ideal such that $A^{\#}=Z$, then $A$ is a
principal ideal of $Q$ and $Q\simeq E_R(R/A)$.

(2) We have $E(U)\simeq Q/I$ for some ideal $I$ of $R$.

(3) We adapt the proof of \cite[Theorem VII.1.4]{fuc}. Suppose that
$\phi:J/I\rightarrow J'/I'$ is an isomorphism. Since $J/I\simeq sJ/sI$
for every $s\in R\setminus Z$ we may assume that $I$ and $I'$ are
proper ideals of $R$. Then $E(J/I)\simeq E(J'/I')$ implies that
$I=(I':t)$ for some $0\ne t\in R$. Since $(J':t)/(I':t)\simeq J'/I'$
we may assume that $I'=I$. By \cite[Theorem VII.1.4]{fuc} we may
assume that $Z\ne0$ and consequently that $R$ is maximal by
\cite[Proposition 1]{gil}.  The isomorphism $\phi$ extends to an
automorphism $\varphi$ of $E(J/I)$. Let $a\in Z$ such that
$Z=(0:a)$. If $I^{\#}=Z$ then $Q=E(J/I)$, and if $I^{\#}\not=Z$ then,
by eventually replacing $I,J$ and $J'$ with respectively
$(I:b),(J:b)$ and $(J':b)$ for some $b\in R$, we may assume that $I$
is faithful and
$Q/Ia=E(J/I)$. By \cite[Corollary VII.2.5]{fuc} $\varphi$ is induced by
multiplication by some $0\ne q\in Q\setminus Z$. Hence $J=(J':q)$. In
the two cases there exists $x\in E(J/I)$ such that
$I=(0:x)=(0:\varphi (x))$. By using Lemma~\ref{L:A} it
follows that $I=(I:q)$.  \qed 
 
\bigskip
From Proposition \ref{P:main} we deduce the following result on the injective
dimension of the $R$-module $R$.
\begin{proposition} Let $R$ be an almost maximal valuation
ring such that $R\ne Q$. Then:
\begin{enumerate}
\item If $Q$ is not coherent then $\mathrm{i.d.}_R\,R=2$.
\item If $Q$ is coherent and not a field then
  $\mathrm{i.d.}_R\,R=\infty$. More precisely, for every $R$-module $M$
  and for every integer $n\geq 1$, then $\mathrm{Ext}^{2n+2}_R(M,R)\simeq
  \mathrm{Ext}^1_R(M,Q/Z)$ and $\mathrm{Ext}^{2n+1}_R(M,R)\simeq\mathrm{Ext}^2_R(M,Q/Z)$.
\end{enumerate}
\end{proposition}
\textbf{Proof.} (1) If we apply Proposition \ref{P:main} to $Q$ we deduce that
$\mathrm{i.d.}_R\,Q/Z=1$. Since $R$ is almost maximal, by Corollary~\ref{C:ifr}
$\mathrm{i.d.}_R\,R/Z=1$. By using the exact sequence: 
$0\rightarrow R/Z\rightarrow Q/Z\rightarrow Q/R\rightarrow
  0$,
we get $\mathrm{i.d.}_R\,Q/R=1$. We deduce the result from this
  following exact sequence:
$0\rightarrow R\rightarrow Q\rightarrow Q/R\rightarrow
  0$.

(2) By using the same exact sequences as in (1) we easily deduce  
that\break 
$\mathrm{Ext}^{p+1}_R(M,R)\simeq \mathrm{Ext}^p_R(M,Q/Z)$, for each $p\geq
3$. On other hand if we apply Proposition~\ref{P:main} to $Q$ we can build
an infinite injective resolution of $Q/Z$ with injective terms
$(E_n)_{n\in\mathbb{N}}$ such that $E_p\simeq H$ if $p$ is even and $E_p\simeq
E(Q/Qa)$, for some $0\not= a\in R$, if $p$ is odd. Now it is easy to
complete the proof. \qed
 
\bigskip
Recall that $\mathrm{i.d.}_R\,R=1$ if and only if $R$ is almost maximal
  when $R$ is a domain which is not a field (\cite[Theorem 4]{mat}) and that
  $\mathrm{i.d.}_R\,R=0$ if and only if $R$ is almost maximal and $R=Q$
  (\cite[Theorem 2.3]{kla}).

Let $U$ be a uniform module over a valuation ring $R$. Recall that if $x$ and
$y$ are nonzero elements of $U$ such that
$(0:x)\subseteq (0:y)$ there exists $t\in R$ such that
$(0:y)=((0:x):t)$: see \cite{nis}. As in \cite[p.144]{fuc}
, we set $U_{\#}=\{s\in R\mid\exists u\in U, u\not= 0 \hbox{ and
  } su=0\}$. Then $U_{\#}$ is a prime ideal and the following lemma holds.
\begin{lemma} 
\label{L:dies}
Let $R$ be a valuation ring and $U$ a uniform
$R$-module. Then for every nonzero element $u$ of $U$,
$U_{\#}=(0:u)^{\#}$.
\end{lemma}
\textbf{Proof.} We set $A=(0:u)$. Let $s\in A^{\#}$. There exists
$t\in (A:s)$ such that $tu\not=0$. We have $stu=0$ whence $s\in
U_{\#}$. Conversely let $s\in U_{\#}$. There exists  $0\not= x\in U$
such that $s\in (0:x)\subseteq
(0:x)^{\#}=A^{\#}$. The last equality holds by Lemma \ref{L:die}.
\qed
 
\begin{proposition} 
\label{P:inj}
Let $R$ be a valuation ring and $U$ a uniform
fp-injective module. Assume that $U_{\#}=Z=P$. Then:
\begin{enumerate}
\item $U$ is faithful when $P$ is finitely generated or faithful.
\item If $P$ is not faithful and not finitely generated then
  $ann(U)=(0:P)$ if $E(U)\not\simeq E(R/P)$ and $U$ is faithful
  if $E(U)\simeq E(R/P)$.
\end{enumerate}
\end{proposition}
\textbf{Proof.} If $E(U)\simeq E(R/P)$ let $u\in U$ such that
$(0:u)=P$. Then for each $0\ne t\in P$, $(0:t)\subseteq
P$. Hence there exists $z\in U$ such that $tz=u\ne 0$. Hence $U$ is
faithful. We assume in the sequel that $E(U)\not\simeq E(R/P)$.

(1) First we suppose that $P=Rp$. Then for every
non-finitely generated ideal $A$ it is easy to check that
$A=(A:p)$. Consequently, for each $u\in U$, $(0:u)$ is
principal, whence $E(U)\simeq E$. Now we assume that $P$ is
faithful. Then $P$ is not principal. For some $0\ne u\in U$ we put
$A=(0:u)$. Let $0\ne t\in P$. Then $(0:t)\subset P$. The
equality $A^{\#}=P$ implies that there exists $s\in P\setminus A$ such
that $(0:t)\subset (A:s)$. We have $su\ne 0$ and
$(0:su)=(A:s)$. It follows that there exists $z\in U$ such
that $tz=su\ne 0$.

(2) We use the same notations as in (1). If
$t\notin (0:P)$ we prove as in the first part of the proof that there
exists $z\in U$ such that $tz\ne 0$. On the other hand for every
$s\notin (0:A)$, $(0:P)\subseteq sA$. Hence $ann(U)=(0:P)$.
\qed
 
\bigskip
From the previous proposition we deduce the following corollary.
\begin{corollary} 
\label{C:inj}
Let $R$ be a valuation ring and $U$ an
indecomposable injective $R$-module. Then the following assertions are
true:
\begin{enumerate}
\item If $Z\subset U_{\#}$ then $U$ is faithful. Moreover, if $R$ is
  almost maximal then $U$ is a factor module of $H$.
\item If $U_{\#}\subset Z$ then $ann(U)=ke(U_{\#})$.
\item If $U_{\#}=Z$ then $U$ is faithful if $Z$ is faithful or
  finitely generated over $Q$. If $Z$ is not faithful and not finitely
  generated over $Q$ then $ann(U)=(0:Z)$ if $U\not\simeq H$ and
  $H$ is faithful.
\end{enumerate}
\end{corollary}

\section{Countably generated indecomposable injective modules}
\label{S:cou}

When $R$ is not a domain we don't know if condition (C) implies
that every indecomposable injective module is countably
generated. However it is possible to give sufficient and necessary
conditions for every indecomposable injective $R$-module to be
countably generated, when $R$ is an almost maximal valuation ring. 

 The following lemmas will be useful in the sequel:
\begin{lemma}
\label{L:cog1} 
Let $R$ be a valuation ring and $A$ a proper ideal
of $R.$ Then $A\not=\bigcap_{r\notin A}Rr$ if and only if there exists
$t\in R$ such that $A=Pt$ and $\bigcap_{r\notin A}Rr=Rt.$
\end{lemma}
\textbf{Proof.} Let $t\in (\bigcap_{r\notin A}Rr)\setminus A.$ Clearly
$Rt=\bigcap_{r\notin A}Rr,$ whence $A=Pt.$\qed
 
Recall that an $R$-module $M$ is \textit{finitely} (respectively
 {\it countably}) {\it cogenerated} if $M$ is a submodule of a product
 of finitely (respectively countably) many injective hulls of simple modules.
\begin{lemma} 
\label{L:cog2}
Let $R$ be a valuation ring and $A$ a proper ideal
of $R$. Suppose that $R/A$ is not finitely cogenerated. Consider the
following conditions:
\begin{enumerate}
\item There exists a countable family $(I_n)_{n\in\mathbb{N}}$ of ideals
  of $R$ such that $A\subset I_{n+1}\subset I_n,\forall n\in\mathbb{N}$
  and $A=\bigcap_{n\in\mathbb{N}}I_n.$
\item There exists a countable family $(a_n)_{n\in\mathbb{N}}$ of elements
  of $R$ such that $A\subset Ra_{n+1}\subset Ra_n,\forall n\in\mathbb{N}$
  and $A=\bigcap_{n\in\mathbb{N}}Ra_n.$ 
\item $R/A$ is countably cogenerated.
\end{enumerate}
Then (1) implies (2) and (3) is equivalent to (2).
\end{lemma}
\textbf{Proof.}  If we take $a_n\in I_n\setminus I_{n+1},\forall
n\in\mathbb{N},$ then $A=\bigcap_{n\in\mathbb{N}}Ra_n.$ Consequently
$(1)\Rightarrow (2)$. It is obvious that $A=\bigcap_{n\in\mathbb{N}}Ra_n$
if and only if $A=\bigcap_{n\in\mathbb{N}}Pa_n,$ and this last
condition is equivalent to: $R/A$ is a submodule of
$\prod_{n\in\mathbb{N}}(R/Pa_n)\subseteq E(R/P)^{\mathbb{N}}$. Hence 
conditions (2) and (3) are equivalent. \qed
 
\begin{lemma} 
\label{L:cog3}
Let $R$ be a ring (not necessarily commutative). Then
the follo\-wing conditions are equivalent.
\begin{enumerate}
\item Every cyclic left $R$-module is countably cogenerated.
\item Each finitely generated left $R$-module is countably
  cogenerated.
\end{enumerate}
\end{lemma}
\textbf{Proof.} Only $(1)\Rightarrow (2)$ requires a proof. Let $M$ be a
left $R$-module generated by $\{x_k\mid 1\leq k\leq p\}$. We induct
on $p$. Let $N$ be the submodule of $M$ generated by $\{x_k\mid 1\leq
k\leq p-1\}$. The induction hypothesis implies that $N$ is a submodule of $G$
and $M/N$ a submodule of $I$, where $G$ and $I$ are product of
countably many injective hulls of simple left $R$-modules. The
inclusion map $N\rightarrow G$ can be extended to a morphism
$\phi:M\rightarrow G$. Let $\varphi$ be the composition map
$M\rightarrow M/N\rightarrow I$. We define $\lambda:M\rightarrow
G\oplus I$ by $\lambda(x)=(\phi(x),\varphi(x))$ for every $x\in M$. It
is easy to prove that $\lambda$ is a monomorphism and conclude the proof.
\qed
 
\begin{proposition} 
\label{P:cog}
Let $R$ be a valuation ring such that
$Z=P$. Consider the following conditions:
\begin{enumerate}
\item $R$ and $R/(0:P)$ are countably cogenerated.
\item $P$ is countably generated.
\item Every indecomposable injective $R$-module $U$ such that
  $U_{\#}=P$ is coun\-tably generated.
\end{enumerate}
Then conditions (1) and (2) are equivalent, and they
are equivalent to (3) when $R$ is almost maximal.

Moreover, when the two first conditions are satisfied, every ideal
$A$ such that $A^{\#}=P$ is countably generated and $R/A$ is
countably cogenerated. 
\end{proposition}
\textbf{Proof.} (1)$\Rightarrow$(2). We may assume that $P$ is not
finitely generated. If $(0:P)=\cap_{n\in\mathbb{N}}Rs_n$, where $s_n\notin
(0:P)$ and $s_n\notin Rs_{n+1}$ for every $n\in\mathbb{N}$, then, by using
 \cite[Proposition 1.3]{kla}, it is easy to prove that
$P=\cup_{n\in\mathbb{N}}(0:s_n)$. Since $(0:s_n)\subset (0:s_{n+1})$ for
each $n\in\mathbb{N}$, we deduce that $P$ is countably generated.

(2)$\Rightarrow$(1). First we assume that $P$ is principal. Then
$(0:P)$ is the nonzero minimal submodule of $R$, and $(0:P^2)/(0:P)$ is
the nonzero minimal submodule of $R/(0:P)$. Hence $R$ and $R/(0:P)$
are finitely cogenerated. Now assume that $P=\cup_{n\in\mathbb{N}}Rt_n$
where $t_{n+1}\notin Rt_n$ for each $n\in\mathbb{N}$. As above we get that
$(0:P)=\cap_{n\in\mathbb{N}}(0:t_n)$. Since $(0:t_{n+1})\subset (0:t_n)$
for each $n\in\mathbb{N}$ it follows that $R/(0:P)$ is countably
cogenerated. If $(0:P)\ne 0$ then $R$ is finitely cogenerated.

(3)$\Rightarrow$(1). It is sufficient to prove that $R/(0:P)$ is
countably cogenerated. We may assume that $P$ is not principal. Then
$F\not\simeq E(R/P)$ and $F_{\#}=P$. Let $\{x_n\mid n\in\mathbb{N}\}$ be a
generating subset of $F$ such that $x_{n+1}\notin Rx_n$ for each
$n\in\mathbb{N}$. By Proposition \ref{P:inj} the following equality holds:
$(0:P)=\cap_{n\in\mathbb{N}}(0:x_n)$. We claim that
$(0:x_{n+1})\subset (0:x_n)$ for each $n\in\mathbb{N}$ else
$Rx_{n+1}=Rx_n$. Consequently $R/(0:P)$ is countably cogenerated.

(1)$\Rightarrow$(3). If $P$ is principal then an ideal $A$ satisfies
$A^{\#}=P$ if and only if $A$ is principal (see the proof of
Proposition \ref{P:inj}). It follows that $U\simeq
R$. Now we suppose that $P$ is not finitely generated. Assume that
there exists $x\in U$ such that
$(0:x)=(0:P)$. If $(0:P)=0$ then $Rx\simeq R$. It follows
that $U=Rx$. If $(0:P)\ne 0$ then $Rx\simeq R/(0:P)$. Since $R$
is not an IF-ring in this case, $R/(0:P)$ is injective by Corollary
\ref{C:ifr}. It follows that $U=Rx$. If $(0:P)\ne 0$ then $E(R/P)\simeq
R$. Hence, if $U$ is not finitely generated, we may assume that
 $(0:P)\subset (0:x)$ for each $x\in U$.  We know that
 $\cap_{n\in\mathbb{N}}Rs_n=(0:P)$ where $s_n\notin (0:P)$ and
 $s_{n+1}\notin Rs_n$ for each $n\in\mathbb{N}$. Let $(x_n)_{n\in\mathbb{N}}$
 a sequence of elements of $U$ obtained by the following way: we pick $x_0$ a
 nonzero element of $U$; by induction on $n$ we pick $x_{n+1}$ such that
 $(0:x_{n+1})\subset (0:x_n)\cap Rs_{n+1}$. This is
 possible since $ann(U)=(0:P)$ by Proposition \ref{P:inj}. Then we
 get that $\cap_{n\in\mathbb{N}}(0:x_n)=(0:P)$. If $x\in U$ then
 there exists $n\in\mathbb{N}$ such that
 $(0:x_n)\subseteq (0:x)$. Hence $x\in Rx_n$ since $U$
 is uniserial.
 
Now we prove the last assertion. If $P$ is principal then $A$ is also principal and $R/A$ finitely
cogenerated. Assume that $P=\cup_{n\in\mathbb{N}}Rs_n$. If $A=Pt$ for some $t\in R$ then $A$ is countably generated and $R/A$
is finitely cogenerated. We may assume that $(A:t)\subset P$ for each
$t\in R\setminus A$. Clearly $A\subseteq\cap_{n\in\mathbb{N}}(A:s_n)$. If
$b\in\cap_{n\in\mathbb{N}}(A:s_n)$, then $b\in (A:P)$ and it follows that
$P\subseteq (A:b)$. Hence $b\in A$, $A=\cap_{n\in\mathbb{N}}(A:s_n)$ and
$R/A$ is countably cogenerated. Let $s\in P\setminus (0:A)$. Thus
$((0:A):s)=(0:sA)\supset (0:A)$. It follows that
$(0:A)^{\#}=P$. Therefore $R/(0:A)$ is countably cogenerated. If
$(0:A)=Pt$ for some $t\in R$, then $tA$ is the nonzero minimal ideal
of $R$ and by using Lemma \ref{L:B} we show that $A$ is principal. If
$(0:A)=\cap_{n\in\mathbb{N}}Rt_n$ then we prove that
$A=\cup_{n\in\mathbb{N}}(0:t_n)$, by using \cite[Proposition 1.3]{kla},
when $A$ is not principal. Hence $A$ is countably generated. \qed
 
\bigskip
Recall that a valuation ring $R$ is \textit{archimedean} if its
maximal ideal $P$ is the only non-zero prime ideal, or equivalently
$\forall a,b\in P,a\ne 0,\exists n\in\mathbb{N}$ such that $b^n\in Ra.$ By
using this last condition we prove that $P$ is countably generated.
\begin{lemma} 
\label{L:arch}
Let $R$ be an archimedean valuation ring. Then its
maximal ideal $P$ is countably generated.
\end{lemma}
\textbf{Proof.} We may assume that $P$ is not finitely generated. Let
$r\in P.$ Then there exist $s$ and $t$ in $P$ such that $r=st$ and
there exists $q\in P$ such that $q\notin Rs\cup Rt$. Hence for each
$r\in P$ there exists $q\in P$ such that $q^2\notin Rr.$ Now we consider the
sequence $(a_n)_{n\in\mathbb{N}}$ of elements of $P$ defined in the
following way: we choose a nonzero element $a_0$ of $P$ and by
induction on $n$ we choose $a_{n+1}$ such that $a_{n+1}^2\notin
Ra_n.$ We deduce that $a_n^{2^n}\notin Ra_0,$ for every integer $n\geq
1.$
Let $b\in P.$ There exists $p\in\mathbb{N}$ such that $b^p\in Ra_0.$ Let
$n$ be an integer such that $2^n\geq p.$ It is easy to check that
$b\in Ra_n.$ Then $\{a_n\mid n\in\mathbb{N}\}$ generates $P.$\qed
 
\bigskip
By using this lemma, we deduce from Proposition \ref{P:cog} the
following corollary. 
\begin{corollary} 
\label{C:cog}
Let $R$ be a valuation ring and
 $N$ its nilradical. Consider the following conditions:
\begin{enumerate}
\item For every prime ideal $J\subseteq Z$, $J$ is countably generated and
 $R/J$ is coun\-tably cogenerated.
\item For every prime ideal $J\subseteq Z$ which is the union of the
  set of primes properly contained in $J$ there is a countable subset
  whose union is $J,$ and for every prime ideal $J\subseteq Z$ which
  is the intersection of the
  set of primes containing properly $J$ there is a countable subset
  whose intersection is $J.$
\item Every indecomposable injective $R$-module is countably generated.
\end{enumerate}
Then conditions (1) and (2) are equivalent and they
are equivalent to (3) when $R$ is almost maximal.

Moreover, when the two first conditions are satisfied, every ideal $A$
of $Q$ is countably generated and $Q/A$ is countably cogenerated.
\end{corollary}
\textbf{Proof.} $(3)\Rightarrow (1).$ For each prime ideal $J\subseteq
Z$, $R_J$ is an indecomposable injective $R$-module. Hence $R_J$ is countably
generated. It is obvious that $J=\cap_{n\in\mathbb{N}}Rt_n$, where
$t_n\notin J$ for each $n\in\mathbb{N}$ if and only if $\{t_n^{-1}\mid
n\in\mathbb{N}\}$ generates $R_J$. Hence $R/J$ is countably
cogenerated. By Proposition \ref{P:cog} $JR_J$ is countably generated over
$R_J$. It follows that $J$ is countably generated over $R$ too.
 
$(1)\Rightarrow (3).$ Since $R/J$ is countably cogenerated and $J$ is
countably generated it follows
that $R_J$ and $JR_J$ are countably generated. By Proposition \ref{P:cog} $U$
is countably generated over $R_J$ and over $R$ too, for every
indecomposable injective $R$-module $U$ such that $U_{\#}=J$. The
result follows from Corollary \ref{C:A}.

$(1)\Rightarrow (2)$. Suppose that $J$ is the union of
the prime ideals properly contained in $J.$ Let $\{a_n\mid n\in\mathbb{N}\}$ be a spanning set of $J$ such that $a_{n+1}\notin Ra_n$ for each
$n\in\mathbb{N}.$ We consider $(I_n)_{n\in\mathbb{N}}$ a sequence of prime
ideals properly contained in $J$ defined in the following way: we
pick $I_0$ such that $a_0\in I_0$ and for every $n\in\mathbb{N}$ we pick
$I_{n+1}$ such that $Ra_{n+1}\cup I_n\subset I_{n+1}.$ Then $J$ is the
union of the family $(I_n)_{n\in\mathbb{N}}$. Now if $J$ is the
intersection of the prime ideals containing properly $J$, in a similar
way we prove that $J$ is the intersection of a countable family of
these prime ideals.

$(2)\Rightarrow (1)$.  By Lemma \ref{L:cog2} we may assume that $V(J)\setminus
  \{J\}$ has a minimal element $I$. If $a\in I\setminus J$ then
  $J=\cap_{n\in\mathbb{N}}Ra^n.$ Now we prove that $J$ is countably
  generated. If $J=N$ then $R_J$ is archimedean. If $J\ne N,$ we may
  assume that $D(J)$ has a maximal element $I$. Then $R_J/IR_J$ is
  archimedean too. In the two cases $JR_J$ is countably generated over
  $R_J$ by Lemma \ref{L:arch}. On the other hand $R/J$ is countably
  cogenerated, whence $R_J$ is countably generated over $R$. Let us
  observe that $JR_J\simeq J/ke(J)$. It follows that $J$ is countably
  generated over $R$ too.

Now we prove the last assertion. We put $J=A^{\#}$. Then $J\subseteq
Z$. By Proposition \ref{P:cog} $A$ is countably generated over
$R_J$. Since $R_J$ is 
countably generated over $R$ it follows that $A$ is countably
generated over $R$ too. On the other hand, since $R_J/AR_J$ is countably
cogenerated, the inclusion $Q/A\subseteq R_J/AR_J$ implies that $Q/A$
is countably cogenerated too, by Lemma \ref{L:cog2}.
\qed
 
\bigskip
From this corollary we deduce the following results:
\begin{corollary} Let $R$ be a valuation ring. Then the following
conditions are equi\-valent:
\begin{enumerate}
\item Every finitely generated $R$-module is countably cogenerated and
  every ideal of $R$ is countably generated.
\item For each prime ideal $J$ which is the union of the
  set of primes properly contained in $J$ there is a countable subset
  whose union is $J,$ and for each prime ideal $J$ which
  is the intersection of the
  set of primes containing properly $J$ there is a countable subset
  whose intersection is $J.$
\end{enumerate}
\end{corollary}
\textbf{Proof.} It is obvious that $(1)\Rightarrow (2)$.

$(2)\Rightarrow (1)$. When $R$ satisfies the condition (D): $Z=P$, this implication holds
by Corollary \ref{C:cog}. Now we return to the general case. Let $A$ be a
non-principal ideal of $R$ and $r\in A,r\ne 0$. Then the factor ring
$R/Rr$ satisfies the condition (D). Hence $A$ is countably generated
and $R/A$ is countably cogenerated. If $R$ is a domain then as in the
proof of Corollary \ref{C:cog} we show that $R$ is countably cogenerated. If
$R$ is not a domain, then $Q$ satisfies (D) and consequently $Q$ is
countably cogenerated over $Q$. By Lemma \ref{L:cog2} $R$ is countably
cogenerated. We conclude by Lemma \ref{L:cog3}.
\qed
 
\begin{corollary} 
\label{C:cnt}
Let $R$ be a valuation ring such that
$\mathrm{Spec}(R)$ is a countable set. Then:
\begin{enumerate}
\item Every finitely generated $R$-module is countably cogenerated.
\item Every ideal is countably generated.
\item Each fp-injective $R$-module is locally injective if and only if
$R$ is almost maximal.
\end{enumerate}
\end{corollary}
\begin{remark} \textnormal{Let $R$ be a valuation ring and $\Gamma(R)$
    its value group. See \cite{shor} for the definition of
    $\Gamma(R)$. If $\mathrm{Spec}(R)$ is countable, then by
    \cite[Theorem 2]{shor} and \cite[Lemma~12.11, p. 243]{sal} we get
    that $\aleph_0\leq\vert\Gamma(R)\vert\leq
    2^{\aleph_0}$. Conversely if $\vert\Gamma(R)\vert\leq\aleph_0$ it
    is obvious that every ideal is countably generated and that each
    finitely generated $R$-module is countably cogenerated.}
\end{remark}

Let us observe that if an almost maximal  valuation ring $R$ satisfies
the conditions of Corollary \ref{C:cog}, then every indecomposable injective
$R$-module $U$ such that $U_{\#}\subset Z$ is flat since $R_J$ is an
IF-ring. It follows that $\mathrm{p.d.}_R\,U\leq 1$ by 
\cite[Proposition 9.8 p.233]{sal}. On the other hand, when $R$ is a
valuation domain that satisfies (C), first it is proved that $\mathrm{p.d.}_R\,Q=1$ and
afterwards, by using \cite[Theorem 2.4 p.76]{fuc}, or by using methods
of R.M. Hamsher in \cite{ham}, it is shown that $Q$ is countably
generated. When $R$ is not a domain it is possible to prove the
following proposition: 
\begin{proposition} Let $R$ be a valuation ring and $J$ a
nonmaximal prime ideal. The following assertions are equivalent:
\begin{enumerate}
\item $R_J$ is countably generated.
\item $\mathrm{p.d.}_R\,R_J=1$.
\end{enumerate}
\end{proposition}
\textbf{Proof.} (1)$\Rightarrow$(2). By \cite[Proposition 9.8 p.233]{sal} ,
$\hbox{p.d.}_R\,R_J=1$ since $R_J$ is flat and countably generated.

(2)$\Rightarrow$(1). If $R'=R/ke(J)$ then $\mathrm{p.d.}_{R'}\,R_J=1$
   since $R_J$ is flat. Then, after eventually replacing $R$ with
   $R'$, we may assume that every element of $R\setminus J$ is not a
   zerodivisor. We use similar methods as in \cite{ham}. If $S$ is a
   multiplicative subset (called semigroup in \cite{ham}) of $R$ then
   $R\setminus S$ is a prime ideal. First, if $s\in P\setminus J$ we
   prove there exists a prime ideal $J'$ such that $s\notin J'$,
   $J\subseteq J'$, $R_{J'}$ is countably generated and
   $\mathrm{p.d.}_R\,(R_J/R_{J'})\leq 1$: we do a similar proof as in
   \cite[Proposition 1.1]{ham}. Now we prove that $R_{J'}/tR_{J'}$ is
   free over $R/tR$, for every non-zerodivisor $t$
   of $R$: we do as in the proof of \cite[Proposition 1.2]{ham}
   . Suppose that $J'\not= J$ and let $t\in J'\setminus
   J$. Since $s$ divides $R_{J'}/tR_{J'}$ and that this module is free
   over $R/tR$, we get that $R_{J'}=tR_{J'}$, whence a contradiction.
\qed

\end{document}